\magnification=1200
\baselineskip=14pt

\def\qed{{$\vrule height4pt depth0pt width4pt$}}
\def\ss{{\smallskip}}
\def\ms{{\medskip}}
\def\bs{{\bigskip}}
\def\ni{{\noindent}}
\def\S{{\Sigma}}

\def\s{{\sigma}}

\def\g{{\gamma}}

\def\l{{\lambda}}

\def\Z{{\bf Z}}

\def\<{{\langle}}
\def\>{{\rangle}}

\def\D{{\Delta}}

\def\R{{\bf R}}
\def\<{{\langle}}
\def\>{{\rangle}}

\def\m1{{\quad ({\rm mod\ 1})}}
\def\L{{\Lambda}}
\def\n{{\bf n}}
\def\Rd{{\cal R}_d}
\def\Sd{{\bf S}^d}
\def\={\ {\buildrel \cdot \over =}\ }

\input epsf
\epsfverbosetrue
 

\centerline{\bf MAHLER MEASURE OF ALEXANDER POLYNOMIALS}
\smallskip
\centerline{Daniel S. Silver and Susan G. Williams} \ss

\ms
\footnote{} {Both authors partially supported by NSF grant
DMS-0071004.}
\footnote{}{2000 {\it Mathematics Subject Classification.}  
Primary 57M25; secondary 37B10, 11R06.}

\noindent {\narrower\narrower\smallskip\noindent  {\bf
ABSTRACT:} Let $l$ be an oriented link of $d$ components in a
homology $3$-sphere. For any nonnegative integer $q$, let $l(q)$
be the link of $d-1$ components obtained from $l$ by performing
$1/q$ surgery on its $d$th component $l_d$. The Mahler
measure of the multivariable Alexander polynomial $\D_{l(q)}$
converges to the Mahler measure of
$\D_l$ as $q$ goes to infinity, provided that
$l_d$ has nonzero linking number with some other
component. If $l_d$ has zero linking number with each
of the other components, then the Mahler measure of $\D_{l(q)}$
has a well defined but different limiting behaviour.  Examples are
given of links $l$ such that the Mahler measure of $\D_l$ is small.
Possible connections with hyperbolic volume are discussed. 


\smallskip}
\ms


\ni {\bf 1. Introduction.} The {\it Mahler measure} of a
nonzero complex Laurent polynomial $f$, introduced by Mahler in
[{\bf Mah60}] and [{\bf Mah62}], is defined by
$${\bf M}(f) = \exp (\int_{\Sd}\log |f({\bf s})|\ d{\bf
s}).$$ Here $d{\bf s}$ indicates integration with respect
to  normalized Haar measure, while $\Sd$ is the 
multiplicative $d$-torus, the subgroup of complex
space
${\bf C} ^d$ consisting of all vectors
${\bf s}=(s_1, \ldots, s_d)$ with $|s_1|= \cdots =|s_d| = 1$.
We adopt the convention that the Mahler measure of the zero polynomial
is $0$.

It is obvious that Mahler measure is multiplicative, and
the measure of any unit is $1$. It is known that $M(f)=1$ if and 
only if $f$ is equal up to a unit factor
to the product of cyclotomic polynomials in a single variable 
evaluated at monomials 
(see [{\bf Sch95}, Lemma 19.1]). 

The quantity ${\bf M}(f)$ is the geometric mean of $|f|$ 
over $\Sd$. By Jensen's formula [{\bf Alh66}, p.
208] the Mahler measure of a nonzero polynomial
$f(u)=c_nu^n + \cdots +c_1u +c_0\
(c_n\ne 0)$ of a single variable is
$${\bf M}(f) = |c_n|\cdot \prod_{j=1}^n {\rm max}(|r_j|,1),$$
where $r_1,\ldots, r_n$ are the roots of $f$. A short proof
can  be found in either [{\bf EW99}] or [{\bf Sch95}]. 

The group ring $\Z[u_1^{\pm 1}, \ldots, u_d^{\pm
1}] \cong \Z\Z^d$ of Laurent polynomials with integer coefficients
will be denoted by $\Rd$. It
is easy to see that the set of Mahler measures of 
polynomials in ${\cal R}_1$ is contained in $\Rd$, for any $d> 1$.
It is widely believed that the containment is proper, but no proof is
known. Furthermore, if $f \in {\cal R}_d$, then there is no obvious
relationship between the Mahler measure of $f$ and that of
the $1$-variable polynomial obtained from $f$ by setting all
of its variables equal.

In
[{\bf Leh33}] Lehmer found that the Mahler measure of  the
degree 10 polynomial
$$L(u)= 1+u-u^3-u^4-u^5-u^6-u^7+u^9+u^{10}$$
is approximately $1.17628$. He subsequently asked whether, 
given any $\epsilon>0$, there exists a polynomial $f \in
{\cal R}_1$ such that $1<M(f)<1+\epsilon$.
Lehmer's question remains open.
Despite extensive computer searches no polynomial $f \in {\cal
R}_1$ has been found such that $1<M(f)<M(L)$. The reader might
consult Waldschmidt [{\bf Wal80}], Boyd [{\bf Boy81}], Stewart
[{\bf Ste78}] or Everest and Ward [{\bf EW99}] for more about
Lehmer's question.

The set ${\cal L} = \{M(f(u)) \mid f \in {\cal R}_1\}$ is a
natural object for investigation. As Boyd 
observed in [{\bf Boy81}], Lehmer's question is equivalent to
the question  of whether  
$1$ is a limit point of ${\cal L}$. The {\it $k$th derived set}
${\cal L}^{(k)}$ is defined inductively: ${\cal L}^{(0)}={\cal
L}$, and ${\cal L}^{(k)}$ is the set of limit points of 
${\cal L}^{(k-1)}$. If $1\in {\cal L}$, then it can be seen from
the multiplicativity of $M$ that ${\cal L}^{(k)} = [1,
\infty),$ for  all $k\ge 1$.

The group ring $\Rd$ has a natural involution \ $\bar{} :\Rd  \to
\Rd$, sending each $u_i$ to $u_i^{-1}$. A polynomial
$f \in {\cal R}_1$ is {\it reciprocal} if $\bar f(u)=
u^n f(u)$ for some $n$. By a
result of Smyth [{\bf Smy71}] the Mahler measure of any
nonreciprocal, irreducible polynomial $f\in {\cal R}_1$ not
divisible by $u-1$ is at least $1.3247\ldots$, the real root of
$u^3-u+1$. The condition that $f$ is not divisible by $u-1$ is
necessary, since multiplying any reciprocal polynomial by $u-1$
yeilds a nonreciprocal polynomial with the same Mahler
measure as $f$. As a consequence of Smyth's theorem,
one need consider only reciprocal polynomials when addressing
Lehmer's question. 

Let $l=l_1 \cup \cdots \cup l_d$	 be 
an oriented link of $d$ components in the $3$-sphere $S^3$ or
more generally an orientable homology $3$-sphere $\S$. (A {\it
homology $3$-sphere} is a closed $3$-manifold with the same
integral homololgy groups as $S^3$.) The {\it exterior}
$E$ of $l$ is the closure of $\S$ minus a tubular neighborhood of
$l$. The homology group $H_1E\cong
\Z^d$ has natural basis represented by the
meridians $m_1, \ldots, m_d$ of $l$ with orientations induced by
the link. All homology groups in this paper have
integer  coefficients.

The abelianization homomorphism $\g:\pi_1E
\to\Z^d$ determines a covering $p:E_\g \to E$, the
universal abelian cover of the link exterior. The ($0th$)
{\it Alexander polynomial} of $l$, denoted here by $\D_l=
\D_l(u_1, \ldots, u_d)$, is the first characteristic polynomial
of the {\it Alexander module}  
$H_1(E_\g, p^{-1}(*))$. (In general, the
$i${\it th Alexander polynomial} is the $(i+1)$th
characteristic polynomial. It is defined only up to
multiplication by units in $\Rd$.) Alexander polynomials are
easily computed from link diagrams. See [{\bf CF63}], [{\bf
Kaw96}].

The following well-known theorem  will be used repeatedly. A
proof can be found in [{\bf Kaw96}]. \bs

\ni {\bf Theorem 1.1. (Torres conditions)} The Alexander polynomial
$\D_l(u_1, \ldots, u_d)$ of a $d$-component link $l=l_1\cup
\ldots\cup l_d$ satisfies: \ms
\item{} (1) $\D_l(u_1,\ldots, u_d)\= \D_l(u_1^{-1}, \ldots,
u_d^{-1})$; \ms
\item{} (2) $\D_l(u_1, \ldots, u_{d-1},1)\=
\cases{\displaystyle{{u_1^{\lambda_1} -1}\over{u_1 -1}}\
\D_{l'}(u_1)& if
$d=2$\cr (u_1^{\lambda_1}\cdots u_{d-1}^{\lambda_{d-1}}-1)\
\D_{l'}(u_1,
\ldots, u_{d-1})& if $d\ge 3.$}$ \bs

\ni Here $l'$ denotes the link $l_1\cup\cdots\cup l_{d-1}$ while 
$\lambda_i$ is the linking number ${\rm Lk}(l_i, l_d)$. The
symbol
$\=$ indicates equality up to a unit factor. \bs 

The first Torres condition states that Alexander knot polynomials
are reciprocal. In view of
Smyth's result noted above, they are a natural source of examples
for the study of Lehmer's question.  Adding to their interest is
the observation of Short and Neumann in [{\bf Kir97}] that
Lehmer's polynomial
$L(u)$ is the Alexander polynomial of a knot $k \subset S^3$,
and hence necessarily of infinitely many knots. One such knot is
the $(-2,3,7)$-pretzel knot [{\bf Hir98}]. 

There are other abelian covers of $E$ besides
the universal abelian one. 
Given a finite-index subgroup $\L \subset \Z^d$ one can consider
the finite-sheeted cover $E_\L$ associated to the homomorphism
$\pi_1E\ {\buildrel \g \over \to}\ \Z^d \to\ \Z^d/\L$, where the 
second map is the natural projection. Such a cover can be
completed to a branched cover $M_\L$. 
In [{\bf SW99}] we proved that for any oriented link $l\subset
S^3$, the Mahler measure of $\D_l$ has a natural
topological interpretation as the exponential rate of growth
of the order of the torsion subgroup $TH_1M_\L$, computed as a
suitably defined measure of $\L$ goes to infinity. (Although
results in [{\bf SW99}] were stated for links in $S^3$, they
generalize easily for links in homology spheres $\S$.) 

In view of what has been said, it is reasonable to expect
that topology---knot theory in particular---can shed some light
on Lehmer's question.
\bigskip

\ni {\bf 2. Surgery on links and Mahler measure limits.} Let $k$
be an oriented knot in a homology $3$-sphere $\S$ with tubular
neighborhood $V=S^1\times D^2$. The homology group $H_1\partial V
\cong \Z^2$ has a natural, well-defined basis represented by an
oriented meridian $m$ and a  longitutde $l$ of $k$. A {\it
longitude} is a simple closed curve in $\partial V$ that is
essential but  null-homologous in $\S- int\ V$, oriented in the
direction of
$k$. 

Let $p,q$ be relatively prime integers. Removing $V$ from $\S$
and reattaching it so that $* \times \partial D^2 \subset
\partial V$ represents $p\cdot [\ m\ ]  + q\cdot [\ l\ ] \in
H_1\partial V$ produces a $3$-manfold $\S' =\S'(k;p/q)$ said to
be {\it obtained from $\S$ by $p/q$ surgery on $k$}. It is well
known that $\S'$ is a homology $3$-sphere if and only if 
$p = \pm 1$. \bs

\ni {\bf Definition 2.1.} Let $l=l_1 \cup \cdots \cup l_d$
be an oriented link in a homology $3$-sphere, and let $q$ be 
a nonzero integer. Then $l(q)$ is the oriented link
$l_1 \cup \cdots \cup l_{d-1}$ regarded in $\S'(l_d; 1/q)$.
\bs

If $l$ is a link in $S^3$ and $l_d$ is unknotted, then one
can describe $l(q)$ simply: Let $D$ be a $2$-disk
that bounds
$l_d$. We can assume that the sublink $l_1 \cup \cdots \cup
l_{d-1}$ intersects $D$ transversely. Then  $l(q)$ is the link
in $S^3$ 
obtained from
$l_1 \cup \cdots \cup l_{d-1}$ by cutting the strands that
pass through $D$, twisting $q$ full times in the direction of
the longitude of $l_d$, and then reconnecting. (Details can
be found in [{\bf Rol76}].) Kidwell investigated this case in
[{\bf Kid82}]. He showed that if $l_d$ has nonzero linking
number with some other component of $l$, then, as $q$ goes to
infinity, the degree of the ``reduced Alexander polynomial''
of $l(q)$, the polynomial obtained from
$\D_{l(q)}$ by setting $u_1 = \cdots = u_{d-1}$, grows
without bound. He found, in fact, that the sequence of nonzero
exponents of the reduced Alexander polynomial acquires
an ever-expanding gap. In contrast, Theorem 2.2
asserts that the Mahler measures of the Alexander polynomials
are well behaved. 
\bs

\ni {\bf Theorem 2.2.} Assume that $l = l_1 \cup \cdots \cup
l_d$ is an oriented link in  a homology $3$-sphere, and let $l(q)$
be as in Definition 2.1. If some  linking number $\l_i = {\rm
Lk}(l_i, l_d)$ is nonzero,
$1\le i\le d-1$, then 
$$\lim_{q\to \infty} M(\D_{l(q)}) = M(\D_l).$$ \bs

\ni {\bf Remarks 2.3.} 
1. By replacing $1/q$ with $p/q$ in Theorem
2.2 more general, albeit more complicated, results are possible.
For the sake of simplicity we have chosen not to work in
such generality. 

2. Theorem 2.2 bears a striking resemblance of form to a
theorem of Thurston [{\bf Thu83}] which states that the volume 
of a hyperbolic $3$-manifold with cusps is the limit of the
volumes of the manifolds obtained by performing $(p_i,q_i)$
Dehn surgery on the $i$th cusp. The pairs $(p_i,q_i)$ are
required to go to infinity in a suitable manner. (See also
[{\bf NZ85}].) We are grateful to John Dean for bringing this to
our attention. \bs

In [{\bf Kid82}] Kidwell did not address the case in which $l_d$
has zero linking number with each of the remaining components of
$l$. Theorem 2.4  completes the picture. \bs

\ni {\bf Theorem 2.4.} Assume that $l=l_1 \cup \cdots
\cup l_d$ is an oriented link in a homology $3$-sphere.
\break 

(1) If
$\l_i=0$,  for each $1\le i \le d-1$, then
$$\lim_{q\to \infty} {1 \over q}\ \D_{l(q)}(u_1,
\cdots, u_{d-1})
\=\ \cases{(u_1-1)\ \displaystyle {{\partial \over \partial
u_2}\biggr\vert_{u_2 =1}\D_l(u_1, u_2)}, & if $d=2$;\cr
\displaystyle{ {\partial \over \partial u_d}\biggr\vert_{u_d=1}\
\D_l(u_1, \ldots, u_d)}, & if $d\ge 3$.}$$
\ni The convergence is in the strong sense that the
polynomials $\D_{l(q)}$ have eventually constant degree, and the 
coefficients of ${1\over q}\D_{l(q)}$ converge to those of the
polynomial on the right. 

(2)\quad If
$(u_d -1)^2$ divides $\D_l(u_1, \ldots, u_d)$, then for every $q,$
$$\D_{l(q)}(u_1, \ldots, u_{d-1})\ =\
\D_{l'}(u_1, \ldots, u_{d-1}),$$ 
where $l' = l_1\cup \cdots \cup l_{d-1}$. \bs

\ni {\bf Remark 2.6.} The hypothesis of Theorem 2.4 that $\l_i=0$ for 
each $i$ implies that $\D_l(u_1, \ldots, u_{d-1}, 1) = 0$
and thus
$u_d-1$ divides
$\D_l$, by the second Torres condition. Consequently, partial
differentiation of $\D_l$ with respect to $u_d$ followed by evaluation
at $u_d =1$ is equivalent to dividing $\D_l$ by $u_d-1$ and then
setting $u_d$ equal to $1$ in the result. In particular, the
operation is well defined even though $\D_l$ is determined only
up to multiplication by a unit in $\Rd$. \bs

By a theorem of Boyd [{\bf Boy98}] the Mahler measure $M(f)$ is a
continuous function of the coefficients of $f$ for polynomials of 
fixed total degree. This yields the following corollary to Theorem
2.4. 
\bs

\ni {\bf Corollary 2.5.} Under the hypothesis of Theorem 2.4, 
$$\lim_{q\to \infty} {1 \over q}\ M(\D_{l(q)})\ =\
M\biggr[\displaystyle{ {\partial \over \partial u_d}\biggr
\vert_{u_d =1}\ \D_l(u_1, \ldots, u_d)}\biggr].$$\bs
Moreover, if $(u_d -1)^2$ divides $\D_l$, then for every $q$,
$$M(\D_{l(q)})\ = \ M(\D_{l'}).$$

\bs

\ni {\bf 3. Proof of Theorem 2.2.} By a result
proved for a special case by Boyd [{\bf Boy81$'$}] and in
general by Lawton [{\bf Law83}], the Mahler measure of any
polynomial
$f \in \Rd$ can be expressed as the limit of Mahler measures of 
polynomials in a single variable. For ${\bf r} = (r_1, \ldots,
r_d)$ and $\n= (n_1, \ldots, n_d) \in \Z^d$ define
$$| \n |=\max \{|n_1|, \ldots, |n_d|\},$$
$$({\bf r},\n) = r_1 n_1+\cdots +r_d n_d ,$$
and
$$\< {\bf r}\> =\min \{ |{\bf n}| : {\bf 0}\ne {\bf n}\in
\Z^d, ( {\bf r},{\bf n} ) =0 \}.$$
Also, let 
$$f_{\bf r}(u) = f(u^{r_1}, \ldots, u^{r_d}).$$\bs

\ni {\bf Lemma 3.1}. [Boyd-Lawton Lemma] For every $f\in
\Rd$,
\quad
$\lim_{ \<{\bf r}\>\to \infty } M(f_{\bf r})=M(f).$\bs

We also need the following consequence of Lemma 3.1. \bs

\ni {\bf Corollary 3.2.} Let $\kappa_1, \ldots,  \kappa_{d-1}$ be
integers, not all zero. For any $f \in \Rd$ and positive integer
$q$, let $f^{(q)}$ be the element of $\R_{d-1}$ defined by 
$$f^{(q)}(u_1, \ldots, u_{d-1}) = f(u_1, \ldots,
u_{d-1},\ (u_1^{ \kappa_1}\cdots u_{d-1}^{ \kappa_{d-1}
})^q).$$ Then $\lim_{q\to \infty}M(f^{(q)}) = M(f).$\bs

\ni {\bf Proof.} By Lemma 3.1, given $\epsilon >0$, there
exists $K>0$ such that 
$|M(f_{\bf s})-M(f)| <\epsilon /2$ whenever ${\bf s}\in
\Z^d$ satisfies $\<{\bf s}\>\ge K$. We will show that for $q\ge
K$ we have
$|M(f^{(q)}) - M(f)| <\epsilon$, from which Corollary 3.2
follows.

Fix $q\ge K$. Choose ${\bf r}\in \Z^{d-1}$ such that $\<{\bf r}\>
\ge 2Kq\cdot \max_i \{\kappa_i\}$, and also such that 
$$|M(f_{\bf r}^{(q)}) - M(f^{(q)})| <\epsilon/2.$$  
Set ${\bf r}^+ = (r_1, \ldots, r_{d-1}, q(\kappa_1 r_1+\cdots +
\kappa_{d-1}r_{d-1})).$ Then $f^{(q)}_{\bf r} = f_{{\bf r}^+}$, so
it suffices to show that
$\<{\bf r}^+\>\ge K$.  Suppose $(\n, {\bf r}^+)=0$, where ${\bf
0}\ne
\n=(n_1,\ldots, n_d)$. Then 
$$\sum_{i=1}^{d-1}(n_i+n_dq\kappa_i)r_i = 0.$$
Hence either (i) $n_i+n_dq\kappa_i =0$, for all $i$; or else
(ii) $|n_j+n_dq\kappa_j|\ge \<{\bf r}\> \ge 2Kq\cdot \max_i
|\kappa_i|$, for some $j$. In case (i) we must have $n_d \ne 0$.
Since we assume that 
$\kappa_k\ne 0$ for some $k$, it follows that $|n_k|\ge q\ge K$ and
so $|\n|\ge K$. In case (ii) either $|n_j|$ or $|n_d|\cdot q \cdot
|\kappa_j|$ is at least $K\cdot q\cdot \max_i |\kappa_i|$, and so
again we have
$|\n| \ge K$. \qed \bs

\ni {\bf Proof of Theorem 2.2.} Our proof closely follows a proof
of the  second Torres condition found in [{\bf Kaw96}]. 

Let $E'$ denote the exterior of $l(q)$ in the homology
$3$-sphere $\S' = \S'(l;l_d, 1/q)$. Let 
$$\g: \pi_1 E' \to 
H_1E' \cong \<u_1, \ldots, u_{d-1}\ :\ [u_i, u_j]=1\ (1\le i<j<
d)\>$$ 
be the  abelianization homomorphism mapping the class
of the
$i$th meridian $m_i$ to $u_i$, $\ 1\le i<d$, and mapping the class
of
$m_d$ to $u_1^{-q\l_1}\cdots u_{d-1}^{-q\l_{d-1}}$. We will denote
the corresponding covering space by $p: E'_\g \to E'$.

The exterior $E$ of the original link $l$ is a subspace of 
$E'$. Let $\nu$ be the natural composite epimorphism
$\pi_1 E \to \pi_1E'\ {\buildrel \g \over \to }\ \Z^{d-1}$. The
total space of the corresponding cover $E_\nu \to E$ can be
identified with the preimage $p^{-1}(E) \subset E'_\g$.

Consider the portion of the homology exact
sequence of the pair $(E'_\g, E_\nu)$:
$$\cdots \to H_2E'_\g {\ \buildrel j_* \over \to\ } H_2(E'_\g,
E_\nu) {\ \buildrel \partial_* \over \to\ } H_1E_\nu
{\ \buildrel i_* \over \to\ } H_1E'_\g {\ \buildrel j_* \over
\to\ } H_1(E'_\g,E_\nu) \to
\cdots\eqno(3.1)$$  Here all homology groups may  be regarded as
${\cal R}_{d-1}$-modules. By the excision isomorphism $H_k(E'_\g,
E_\nu)$ is trivial unless $k=2$, and
$H_2(E'_\g, E_\nu)\cong {\cal R}_{d-1}/(u_1^{\l_1}\cdots
u_{d-1}^{\l_{d-1}}-1)$. Hence we have a short exact sequence
$$0 \to ker\ i_* \to H_1E_\nu \to H_1E'_\g \to 0.$$
The $0th$ characteristic polynomials satisfy

$$\D_0(H_1E_\nu)\ \=\ \D_0
(ker\ i_*)\
\D_0(H_1E'_\g).\eqno{(3.2)}$$
As in [{\bf Kaw96}, Proposition 7.3.10], we find that
$\D_0(H_1E'_\g)$ is the Alexander polynomial of $l(q)$, while
$$\D_0(H_1E_\nu)\ \=\ \cases{\D_l^{(q)}&if
$d\ge 3$ \cr
(u_1-1)\ \D_l^{(q)}& if
$d=2$.}$$
(Here we use the
notation of Corollary 3.2 with
$\kappa_i = -\lambda_i,\ 1\le i<d$.) Since
$\D_0(ker\ i_*)$ is a divisor of $\D_0(H_2(E'_\g, E_\nu))
\ {\buildrel \cdot \over =\ } u_1^{\l_1}\cdots
u_{d-1}^{\l_{d-1}}-1,$ its Mahler measure is $1$.
Hence the Mahler measure of  
$\D_{l(q)}$ is equal to that of $\D_l^{(q)}$. Theorem 2.2
now follows from  Corollary 3.2. \qed \bs

\ni {\bf Remark 3.3.} The equation (3.2) can be improved. 
We have
$$\D_l^{(q)}(u_1,\ldots, u_{d-1})\ \=\
\cases{( u_1^{\l_1}\cdots u_{d-1}^{\l_{d-1}}-1)\ \D_{l(q)}&if
$d\ge 3$ \cr
\displaystyle{{u_1^{\l_1}-1} \over 
{u_1-1}}\ 
\D_{l(q)}& if
$d=2$.}\eqno(3.3)$$
If $H_2(E_\gamma')=0$, then (3.3)
follows immediately from (3.1) and (3.2). If $H_2(E_\gamma')\ne 0$,
then as in [{\bf Kaw96}, 7.3.5] we have
$\D_0(H_1E'_\g)=0$; in this case both sides of (3.2) vanish, and
(3.3) is trivial. 

Equation (3.3) can also be obtained by applying Theorem 6.7
of [{\bf Fox60}].
 
\bs

\ni {\bf 4. Proof of Theorem 2.4.} Following [{\bf Kid82}] we
add an unknotted oriented component $l_{d+1}$ to $l=l_1\cup \cdots
\cup l_d$  such that $Lk(l_d, l_{d+1})=Lk(l_{d-1}, l_{d+1})=1$, while
${\rm Lk}(d_i, d_{i+1})=0$ for $1\le i < d-1$.   We
denote the augmented link by $l^+$. By the second Torres condition,
$\D_{l^+}(u_1, \ldots, u_d, 1)\ \=\ (u_{d-1}u_d-1)\
\D_l(u_1, \ldots, u_d).$ Differentiating each side of the equation
and recalling that $\D_l(u_1, \ldots ,u_d)=0$, we have
$${\partial \over \partial u_d}\biggr\vert_{u_d=1}\D_{l^+}(u_1,
\ldots, u_d, 1)\ \=\ (u_{d-1}-1)\ {\partial
\over \partial u_d}\biggr \vert_{u_d=1} \D_l (u_1, \ldots,
u_d).\eqno (4.1)$$

Equation (3.3) implies that
$$\D_{l^+}(u_1, \ldots, u_{d-1}, u_{d+1}^q, u_{d+1})\ \=\
(u_{d+1}-1)\ \D_{l^+(q)}(u_1, \ldots, u_{d-1}, u_{d+1}),$$
where $l^+(q)$ is the $d$-component link obtained from  $l^+$ by
performing $1/q$-surgery on the component $l_d$.
Again by differentiating, and applying the second Torres condition,
we have
$${\partial \over \partial
u_{d+1}}\biggr\vert_{u_{d+1}=1}\D_{l^+}(u_1, \ldots, u_{d-1},
u_{d+1}^q, u_{d+1})\ \=\ \D_{l^+(q)}(u_1,
\ldots, u_{d-1}, 1)$$
$$\=\ \cases{\D_{l(q)}(u_1)&if $d=2$\cr
(u_{d-1}-1)\ \D_{l(q)}(u_1,\ldots, u_{d-1}) &if $d \ge
3$.}\eqno(4.2)$$

Comparing (4.1) and (4.2) we see that in order to prove the
first assertion of Theorem 2.4 it suffices to show that
$$\lim_{q\to\infty} {1\over q}\ {\partial \over \partial
u_{d+1}}\biggr\vert_{u_{d+1}=1}\D_{l^+}(u_1, \ldots, u_{d-1},
u_{d+1}^q, u_{d+1})\ \=\ {\partial \over
\partial u_d}\biggr\vert_{u_d=1}\D_{l^+}(u_1, \ldots, u_d,
1).\eqno(4.3)$$
By collecting terms, we can write the Alexander polynomial of
$\D_{l^+}$ in the form
$$\D_{l^+}(u_1, \ldots, u_{d+1})\ \=\
\sum_{i=0}^{m}\sum_{j=0}^n f_{ij}(u_1, \ldots, u_{d-1})
\ u_d^i \ u_{d+1}^j$$
for suitable $f_{ij}\in {\cal R}_{d-1}$. A simple calculation
shows that $${\partial \over \partial
u_{d+1}}\biggr\vert_{u_{d+1}=1}\D_{l^+}(u_1, \ldots, u_{d-1},
u_{d+1}^q, u_{d+1})\ \=\sum_{i=0}^m\sum_{j=0}^n
f_{ij}\cdot(qi+j)$$ and
$${\partial \over
\partial u_d}\biggr\vert_{u_d=1}\D_{l^+}(u_1, \ldots, u_d,
1)\ \= \sum_{i=0}^m\sum_{j=0}^n f_{ij}\cdot i,$$
so that (4.3) is immediate.

Now suppose that
$(u_d-1)^2$ divides $\D_l$. Then from (4.1),
$$\sum_{i=0}^m\sum_{j=0}^nf_{ij}\cdot i\ \=\ {\partial\over \partial
u_d}\biggr\vert_{u_d=1} \D_{l^+}(u_1,
\ldots,u_d, 1)=0.$$
Thus equation (4.3) bcomes
$$\sum_{i=0}^m\sum_{j=0}^nf_{ij}\cdot j\ \=\
\cases{\D_{l(q)}(u_1)&if
$d=2$\cr (u_{d-1}-1)\ \D_{l(q)}(u_1,\ldots, u_{d-1}) &if $d \ge
3$.}$$

Let $\tilde l\  = \ l_1\cup\cdots\cup l_{d-1}\cup
l_{d+1}$. By the second Torres condition:
$$\D_{l^+}(u_1, \ldots, u_{d-1}, 1, u_{d+1})\ \=\ (u_{d+1}-1)\
\D_{\tilde l}(u_1, \ldots, u_{d-1}, u_{d+1}).$$
Differentiating, we have
$$\D_{\tilde l}(u_1, \ldots, u_{d-1}, 1)\ \=\ {\partial \over
\partial u_{d+1}}\biggr\vert_{u_{d+1}=1}\D_{l^+}(u_1, \ldots,
u_{d-1}, 1, u_{d+1})\ \=\ \sum_{i=1}^m\sum_{j=1}^n
f_{ij}\cdot j.$$
Hence
$$\D_{\tilde l}(u_1, \ldots, u_{d-1}, 1)\ \=\
\cases{\D_{l(q)}(u_1)&if $d=2$\cr (u_{d-1}-1)\D_{l(q)}(u_1,
\ldots, u_{d-1})&if $d\ge 3$.}\eqno(4.4)$$
On the other hand,
$$\D_{\tilde l}(u_1, \ldots, u_{d-1}, 1)\ \=\
\cases{\D_{l'}(u_1)&if $d=2$\cr (u_{d-1}-1)\D_{l'}(u_1,
\ldots, u_{d-1})&if $d\ge 3$,}\eqno(4.5)$$
again using the second Torres condition. Comparing (4.4) and (4.5)
we are done.\qed \bs

Examples illustrating Theorem 2.2 can be found in section 5. We
conclude this section with two examples that illustrate Theorem
2.4 and its corollary.

The following lemma of Boyd [{\bf Boy81}] is often useful when
computing Mahler measures of Alexander polynomials. A proof
can also be found in [{\bf Sch95}, p.157]. \bs

\ni {\bf Lemma 4.1.} If $g \in {\cal
R}_{d+1}$ is defined by 
$$g(u_1, \ldots, u_d, u_{d+1}) = u_{d+1}f(u_1, \ldots,
u_d)+f(u_1^{-1}, \ldots, u_d^{-1}),$$
for $f\in  \Rd$, then $M(g) = M(f)$. 
\bs

\ni {\bf Example 4.2.} Consider the Whitehead link $l=5^2_1$ in
Figure 1. Its Alexander polynomial is  $(u_1-1)(u_2-1)$. The knots
$l(q)$ are often called ``twist knots,'' and their Alexander
polynomials $qu_1^2-(2q+1)u_1+q$ are easily
computed.  The Mahler measure of $\D_{l(q)}$
increases without bound as $q$ goes to infinity. However,
$${1\over q}M(\D_{l(q)}) = M\biggr(u_1^2-{{2q+1}\over q} \ u_1
+1\biggr )$$ approaches $M(u_1^2-u_1+1)=1=M(\D_l)$, as predicted by
Corollary 2.5.\bs

\epsfxsize=1.2truein  
\centerline{\epsfbox{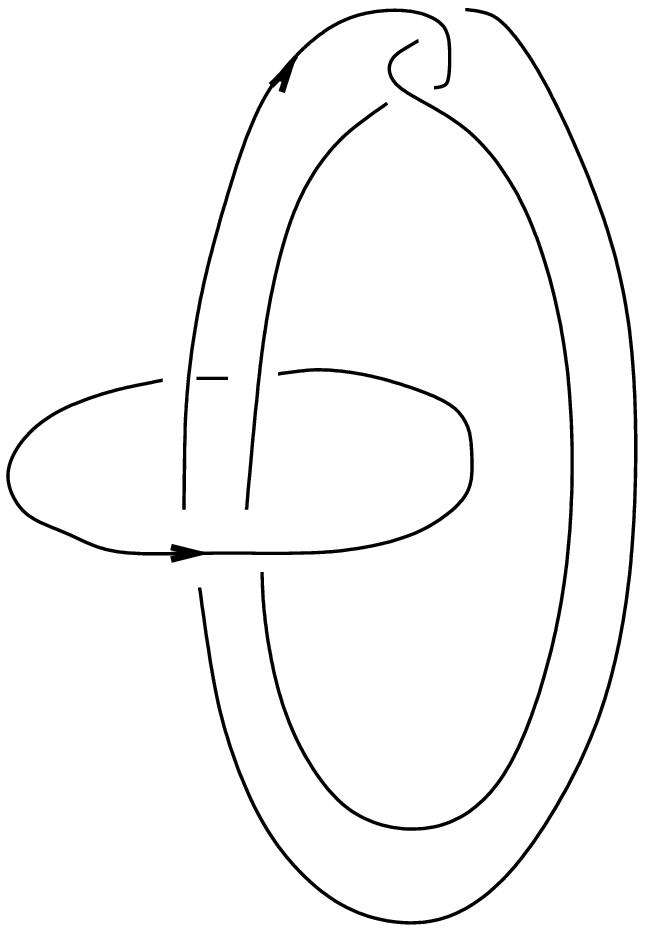}}\bs

\centerline{{\bf Figure 1:} Whitehead link $l=5^2_1$}\bs

\ni {\bf Example 4.3.} Consider the $3$-component link $9^3_8$
in Figure 2. Its Alexander polynomial is
$(u_2-1)(u_3-1)(u_1+2u_2-2u_1u_2-u_2^2).$
Differentiating, we find
$${\partial\over \partial u_3}\biggr\vert_{u_3=1}\D_l(u_1, u_2,
u_3)\= (u_2-1)(u_1+2u_2-2u_1u_2-u_2^2),$$
which can be rewritten as
$(u_2-1)[2u_2-u_2^2-u_1u_2^2(2u_2^{-1}-u_2^{-2})].$
Using Lemma 4.1 and a change of variable we see that 
$$M\biggr[{\partial\over \partial u_3}\biggr\vert_{u_3=1}\D_l(u_1,
u_2, u_3)\biggr]= M[2u_2-u_2^2]=2.$$
By Corollary 2.5 the Mahler measure of $\D_{l(q)}$ is asymptotic to
$2q$. \bs

\epsfxsize=2truein  
\centerline{\epsfbox{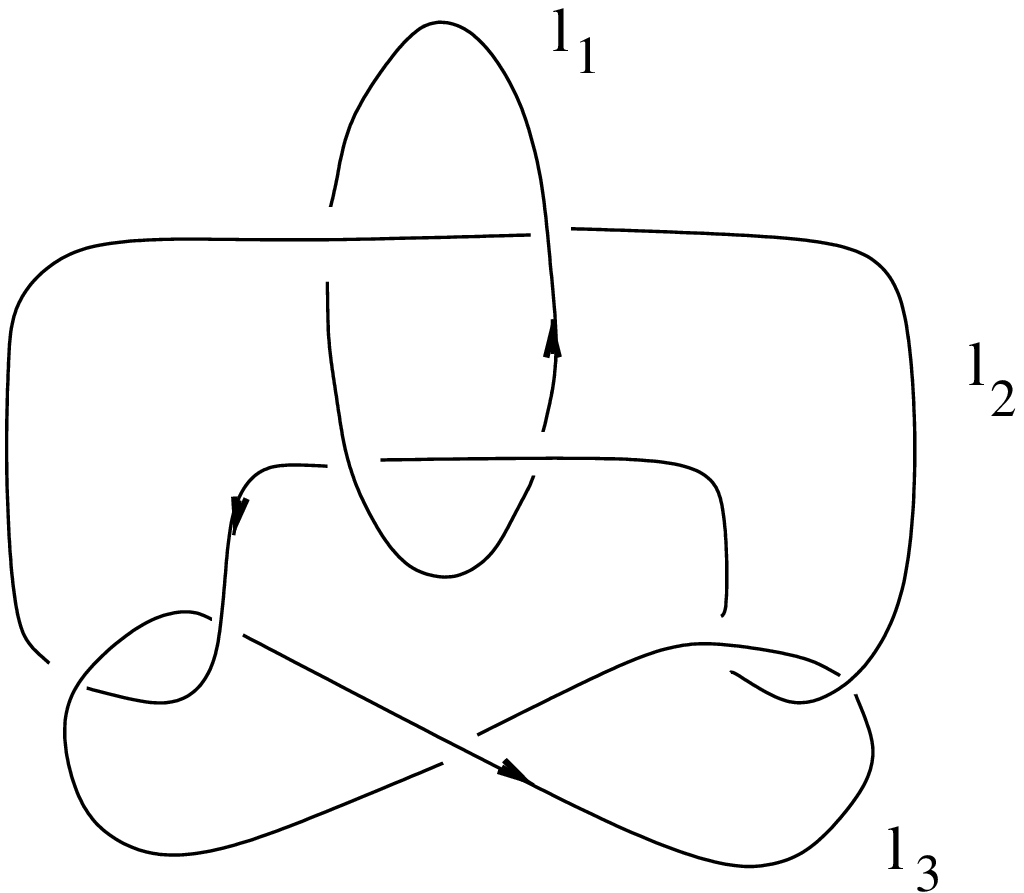}}\bs

\centerline{{\bf Figure 2:} The link $l=9_8^3$}\bs

Examples of $d$-component links $l$ that have Alexander
polynomials divisible by
$(u_d-1)^2$ are not difficult to find. Perhaps the simplest
example is $l=8^2_{10}$. (Here the component $l_1$ is unknotted.) The
reader interested in an exercise can  draw the knots $l(q)$ and
verify that their Alexander polynomials are trivial. \bs

\ni {\bf 5. Alexander link polynomials with small Mahler
measure.}
\bs

\ni {\bf Example 5.1.} Consider the $2$-component link
$l=7^2_1$ in Figure 3a. Its Alexander polynomial $1-u_1 +
(-1+u_1-u_1^2)u_2 + (-u_1+u_1^2)u_2^2$ can be represented
schematically: 
$$\matrix {\hfil & -1 & +1 \cr -1 & +1 &-1\cr +1 &-1 &\cr}$$
Here the number in the $(i+1)$th column from the left and
the $(j+1)$th row from the bottom is the coefficient of 
$u_1^iu_2^j$. 

Replacing 
$u_i$ by $-u_i$, for $i=1,2$, a change that leaves the Mahler
measure unaffected, produces $1+u_1 + (1 + u_1 + u_1^2)u_2 +
(u_1 + u_1^2) u_2$:
$$\matrix {\hfil & +1 & +1 \cr +1 & +1 &+1\cr +1 &+1 &\cr}.$$
Boyd [{\bf Boy78}] has computed the Mahler measure as
approximately $1.25543$.  It is the smallest known
value the derived set of ${\cal L}^{(1)}$ (see [{\bf
Boy78}].) 

The link $l$ is redrawn in Figure 3b so that the second
component $l_2$ appears as a standard unknotted circle. 
The knots $l(q)\subset S^3$ are now easy to visualize: they
are obtained from $l_1$ by giving $q$ full right-hand twists
to the strands passing through $l_2$. By the proof of Theorem
2.2, $\D_{l(q)}$ has the same Mahler measure as $\D_l(u, u^q)$.
When $q =11$, we obtain Lehmer's value $M(L(u))$. When $q=10$
we get the value $1.18836\ldots$, the second smallest known
value of ${\cal L}$ that is greater than $1$. \bs

\epsfxsize=2.3truein  
\centerline{\epsfbox{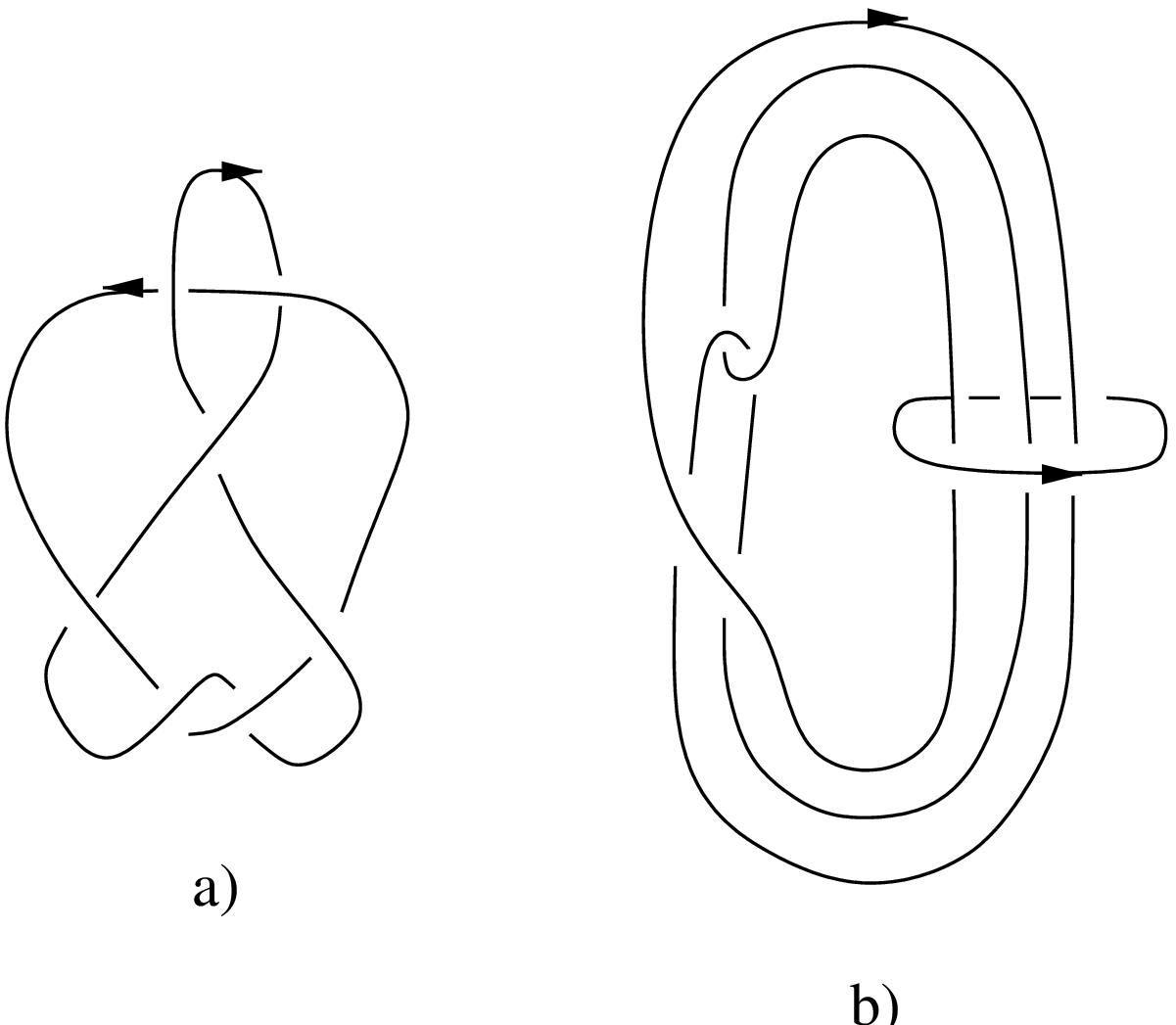}}\bs

\centerline{{\bf Figure 3:} The link $7_1^2$}\bs

\ni {\bf Example 5.2.} The link $6_2^2$ is shown in Figure 2.
Its Alexander polynomial $u_1+(1-u_1+u_1^2)u_2 +u_1u_2^2$:
$$\matrix{&+1&\cr +1&-1&+1 \cr &+1&\cr}$$
has Mahler measure $1.28573\ldots,$ which is the second
smallest known value of ${\cal L}^{(1)}$  (see [{\bf Boy78}]). 
The links $9_4^2,9_9^2$ and $9^2_{50}$ also have Alexander
polynomials with this Mahler measure. 
\bs
\epsfxsize=1.2truein  
\centerline{\epsfbox{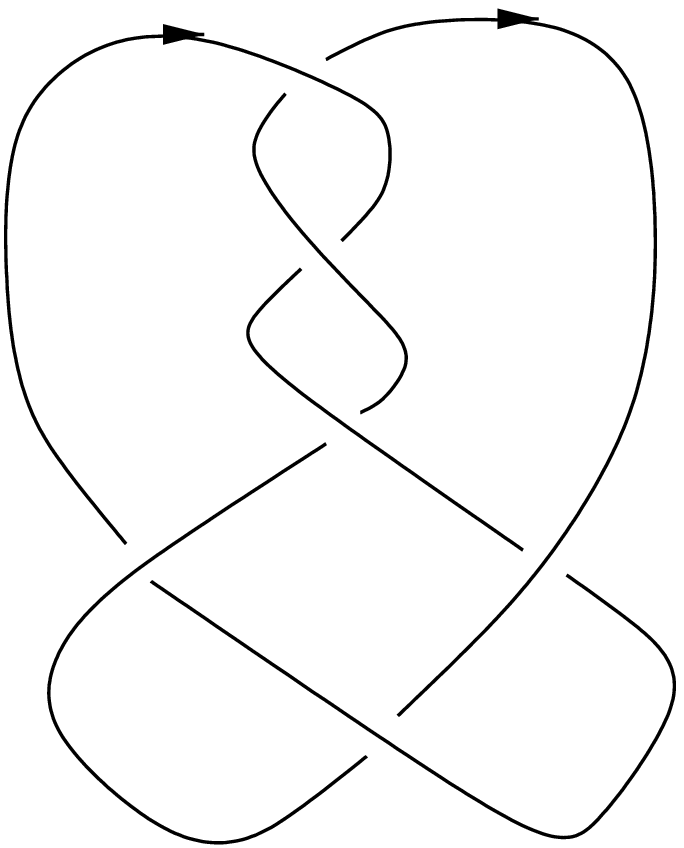}}\bs

\centerline{{\bf Figure 4:} The link $6_2^2$} \bs

\ni {\bf Example 5.3.} Using a computer search, Mossinghoff [{\bf
Mos98}] discovered the polynomial $u_1^2 u_2(u_2+1) +
u_1(u_2^4-u_2^2 +1) + u_2^2(u_2+1)$:
$$\matrix{& +1 & \cr +1 & &&\cr +1& -1& +1\cr&&+1\cr & +1&&\cr}$$
which has Mahler measure approximately
$1.30909$. This value is the third smallest known value
of ${\cal L}^{(1)}$. In a recent private communication
Mossinghoff observed that the symmetric polynomial $$\D(u_1,
u_2)\= u_1^{-2}u_2^{-2}+u_1^{-1}-u_2^{-1}-1+u_1 -u_2+u_1^2u_2^2,$$
schematically:
$$\matrix{&&&&+1\cr&&-1&&\cr&+1&-1&+1&\cr&&-1&&\cr+1&&&&\cr}$$
has the same small Mahler measure. By a
theorem of Levine [{\bf Lev67}] $\D$ is the Alexander polynomial
of a (nonunique)
$2$-component link $l$ in the $3$-sphere. \bs

\ni {\bf Example 5.4.} Figure 5 displays the $3$-component link
$l=6_1^3$. Its Alexander polynomial $\D_l$ is $u_1 + u_2 +
u_3-u_1u_2-u_1u_3-u_2u_3$.
We can rewrite the polynomial as
$u_1 + u_2 - u_1u_2 - u_1u_2u_3 (u_1^{-1}+u_2^{-1} -
u_1^{-1}u_2^{-1})$. We replace
$-u_1u_2u_3$ by $u_3$, a change of variables that does not
affect Mahler measure. 
Lemma 4.1 implies that $M(\D_l)$ is equal to the Mahler measure
of $u_1^{-1} + u_2^{-1} - u_1^{-1} u_2^{-1}$.
Now replacing each $u_i$ by $-u_i$ and multiplying by
the unit $-u_1u_2$ produces the relatively simple polynomial $1 +
u_1 + u_2$ with the same Mahler measure  as $\D_l$. 

\epsfxsize=1.7truein  
\centerline{\epsfbox{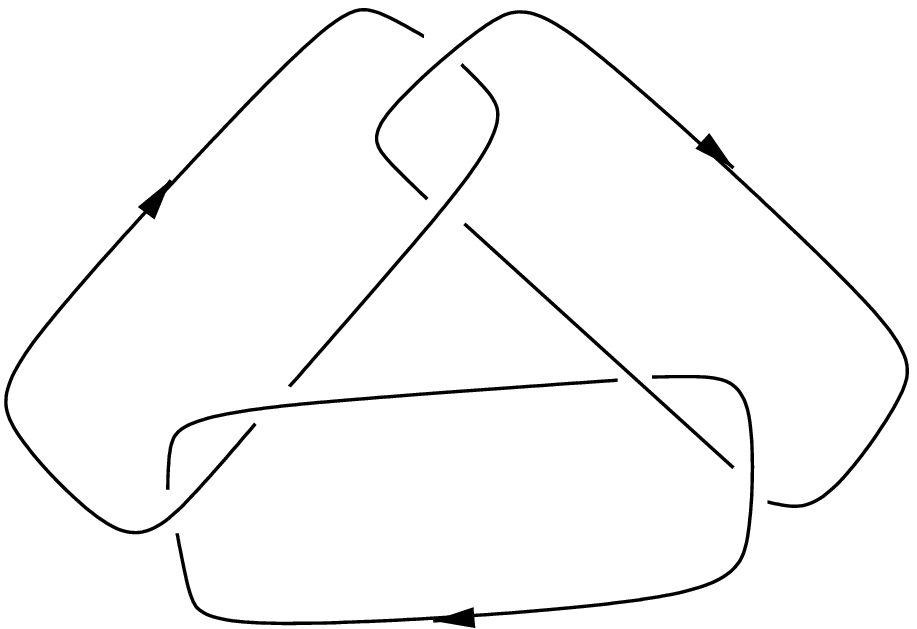}}\bs

\centerline{{\bf Figure 5:} The link $l=6_1^3$.}\bs

Smyth (see [{\bf Boy81}])  has shown that the Mahler
measure of 
$1+ u_1 + u_2$ is 
$$\exp{3 \sqrt 3 \over 4\pi} \sum_{n=1}^\infty {\chi(n) \over
n^2},$$ where $\chi(n)$ is the Legendre symbol
$$\chi (n) = \cases{\hskip .25 cm 1 &if $n\equiv 1\ {\rm mod}$\
3;\cr -1 &if $n\equiv 2\ {\rm mod}$\ 3;\cr
\hskip .3cm 0 &if $n\equiv 0\ {\rm mod}$\ 3\ .\cr}$$
(See also [{\bf EW99}, p. 55] for the calculation.) The Mahler
measures of many polynomials turn out to be simple multiples of
$L$-series of this sort. Deninger [{\bf Den97}] offers
a conjectural explanation in terms of K-theory.

The value of $M(\D_l)$, approximately $1.38135$, is the smallest
known element of ${\cal L}^{(2)}$. Applying Theorem
2.2. twice one can approximate it by Mahler measures of
Alexander polynomials of knots in $S^3$. 
Beginning with $l$, we twist $q_1$ times about $l_3$, producing
the $2$-component link $l(q_1) = l_1' \cup
l_2'$. 
Notice that the components of  $l_1', l_2'$ are each unknotted
and they have nonzero linking number (equal to $+1$). Now we 
twist $q_2$ times about $l_2'$, obtaining  a knot
$k= k(q_1,q_2)$. By Theorem 2.2 the limit
$$\lim_{q_1 \to \infty}\lim_{q_2 \to \infty} M(\D_k)$$
is equal to the Mahler measure $M(\D_l)$ of the link $l=6_1^3$.
The knot $k$  is the closure ot the $(q_1+1)$-braid
$\s_1^{-1}\s_2\cdots \s_{q_1} c^{q_2}$, where $\s_1, \ldots,
\s_{q_1}$ are the usual braid generators, and $c$ is a
full right-hand twist $(\s_1\cdots \s_{q_1})^{q_1+1}.$\bs

\ni {\bf Example 5.5.} Consider the $4$-component link $l= 8^4_2$
in Figure 6. Its Alexander polynomial can be put into  the form
$1-u_1-u_2+u_2u_3+u_1u_2u_3u_4^{-1}(1-u_1^{-1}-u_2^{-1}+
u_2^{-1}u_3^{-1})$.
By Lemma 4.4 the Mahler measure of $\D_l$ is equal to that of
$1-u_1-u_2+u_2u_3$. We replace $-u_1,
-u_2, u_2u_3$ by $u_1,u_2,u_3$, respectively, a change of
variables that does not affect Mahler measure. We then find that
the Mahler measure of
$\D_l$ is equal to that of $1+u_1+u_2+u_3$. By a result of
Smyth (reported in [{\bf Boy81}]) this value is precisely $$\exp
{7 \over 2 \pi^2} \sum_{n=1}^\infty {1\over n^3}.$$ Until
recently,  this was the only nontrivial Mahler measure of a
3-variable polynomial that was evaluated in closed form (see
[{\bf Smy00}]). Its numerical value is  approximately $1.53154$.
\bs
\epsfxsize=1.8truein  
\centerline{\epsfbox{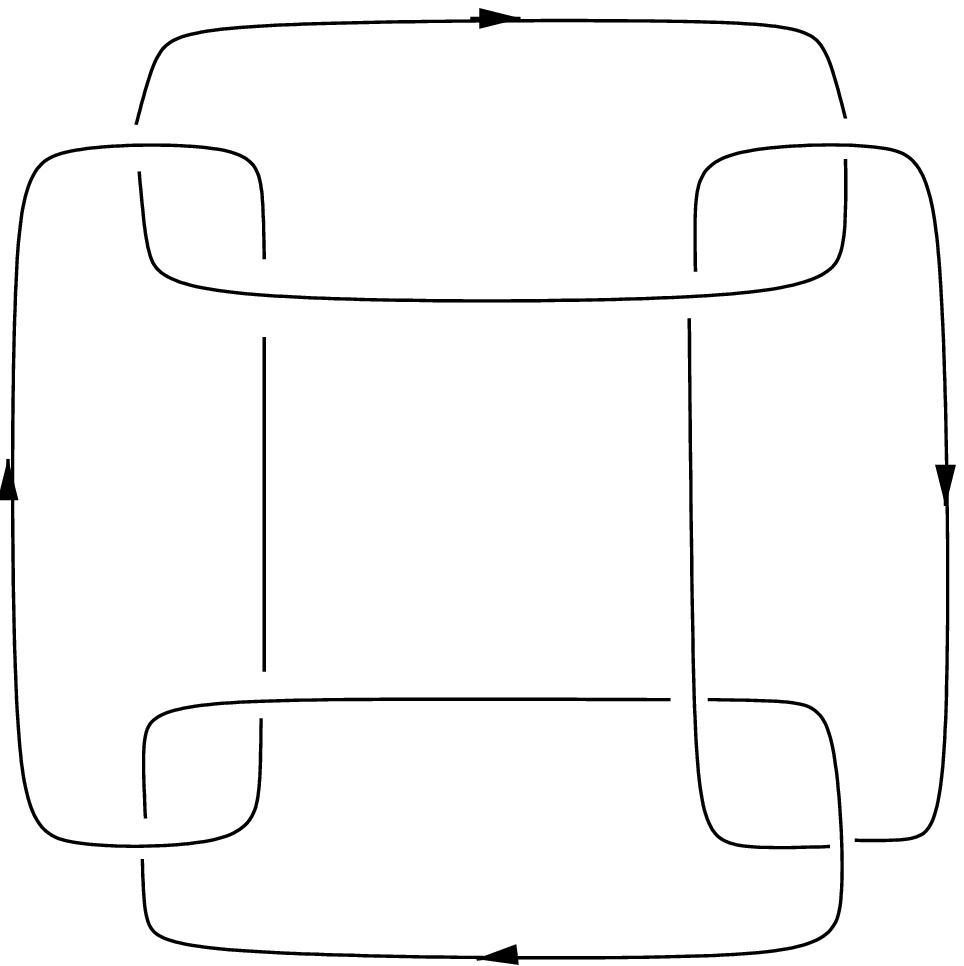}}\bs

\centerline{{\bf Figure 6:}  The link $l=8^4_2$.}
\bs

The links $6_1^3$ and $8_2^4$ are 
examples of pretzel links. For any integers
$p_1. \ldots, p_d$, the {\it pretzel link} $l(p_1, \ldots,
p_d)$ is the boundary of the surface consisting of
two  disks joined by $d$ twisted vertical bands, as in Figure 5.
The $i$th band contains $|p_i|$ half-twists, right-handed if
$p_i$ is positive, left-handed otherwise. The links
$6_1^3$ and $8_2^4$ are $l(2,2,2)$ and $l(2,2,2,-2)$,
respectively. 

The Alexander polynomial $\D_l$ of the pretzel link
$l=l(2,-2,2,-2,2)$ (see Figure 5) has the form
$f + v \bar f,$ where $f(u_1, u_2, u_3, u_4) =
u_1-u_1u_3-u_1u_4-u_2u_3+u_3u_4+u_1u_2u_3$ and $v=
-u_1u_2u_3u_4u_5$. By Lemma 4.4 the Mahler measure of 
$\D_l$ is equal to that of $f$.
Dividing $f$ by $u_1$ produces
$1 -u_3-u_4+u_2u_3+u_1^{-1}u_3u_4-u_1^{-1}u_2u_3$. A further
substitution, replacing $-u_3, -u_4, u_2u_3, u_1^{-1}u_3u_4$ by
$v_1,v_2,v_3,v_4$, respectively, yields
$1+v_1+v_2+v_3+v_4-v_1^{-1}v_2^{-1}v_3v_4$, which has the same
Mahler measure as $f$ and hence as $\D_l$. 

\bs
\epsfxsize=2.4truein  
\centerline{\epsfbox{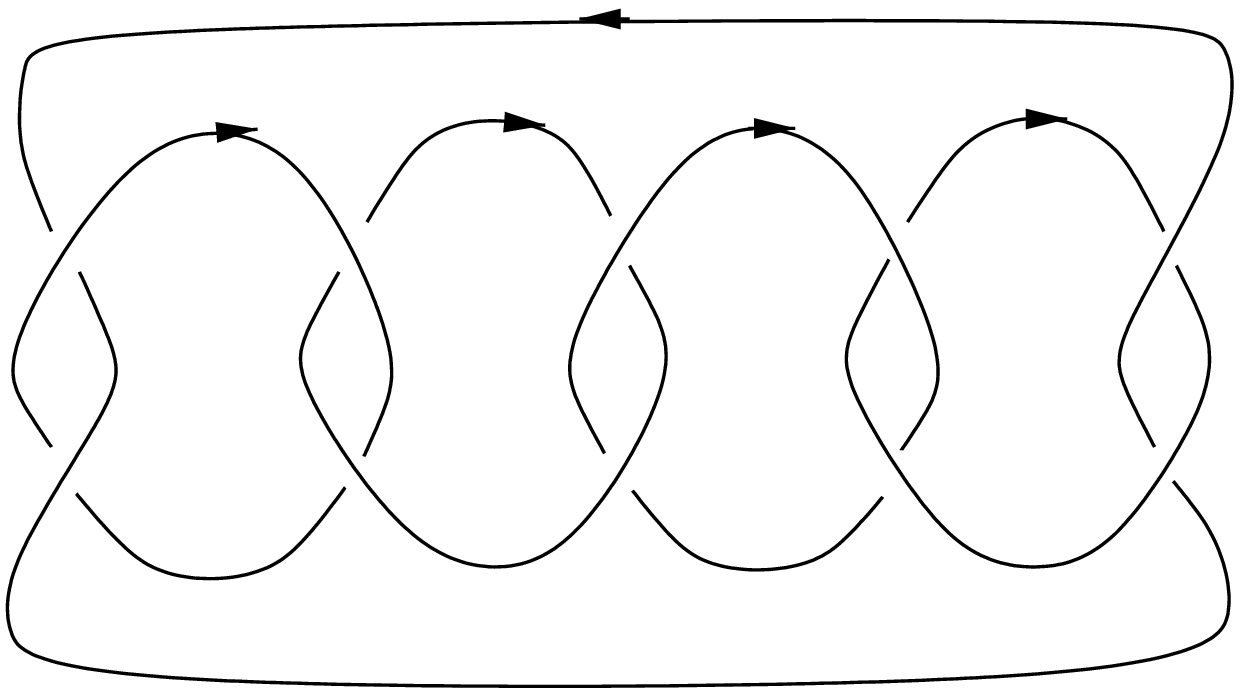}}\bs

\centerline{{\bf Figure 7:} The pretzel link $l(2,-2,2,-2,2)$}\bs

The values of $M(1+ u_1 +u_2 + u_3 + u_4)$ and
$M(\D_l)=M(1+v_1+v_2+v_3+v_4-v_1^{-1}v_2^{-1}v_3v_4)$ are
very close: our calculations suggest that the first
is $1.723\ldots,$ while the second is approximately $1.729$.

The link $l=l(2,-2,2,-2,2)$ has arisen recently in investigations
of  C.~C. Adams [{\bf Ada00}] as well as the work of  N.~M.
Dunfield  and W.~P. Thurston [{\bf Dun00}].  Adams has shown that
$S^3-l$ has the largest cusp density possible for a cusped
hyperbolic $3$-manifold. P.~J. Callahan, C.~D. Hodgson and J.~R.
Weeks have observed that $S^3-l$ is obtained from most of the
small-volume hyperbolic $3$-manifolds of the census [{\bf
HW89}] by removing a shortest-length geodesic and then repeating
the process four additional times. 

Smyth and Myerson [{\bf SM82}] have shown that $M(1 + u_1 +
u_2 + \cdots + u_n)$ is asymptotic to $e^{ -C_0 \sqrt n}$,
where $C_0$ is Euler's constant. Boyd [{\bf Boy00}] has asked
whether this value  is the  rate of growth of $\min {\cal
L}^{(n)}$. \bs

\ni {\bf Question 5.6.} Does there exist a sequence $l_n$ of
$n$-component links such that the sequence of Mahler
measures$M(\D_{l_n})$ is asymptotic to $\min {\cal L}^{(n)}$? \bs

\ni {\bf Definition 5.7.} A number $\theta >1$ is a 
{\it Pisot-Vijayaraghavan number}, or simply a {\it PV
number}, if it is the root of a monic irreducible polynomial
with integer coefficients such that all of its other roots
lie strictly inside the unit disk. If some root lies on the
circle but no other root is outside, then $\theta$ is a {\it
Salem number}. \bs

In [{\bf Sal44}] R. Salem proved that the set of PV numbers is
closed, and hence nowhere dense as it is countable. Since
$1$ is not a PV number, it follows at once that there is a minimum 
PV number $\theta_0$. The value 
of $\theta_0$ was shown by C.~L. Siegel in [{\bf Sie44}] to 
be the real root of $u^3-u-1$, approximately $1.32471$.  

Hironaka proved that among the
reduced Alexander polynomials of pretzel knots and links
$l(p_1,p_2,\cdots, p_k, -1, \cdots, -1)$, where 
$p_1,\ldots, p_k$ are positive and $-1$ occurs $k-2$ times,
Mahler measure is minimized by the Lehmer polynomial $L\=
\D_{l(-2,3,7)}$. (Recall that the  reduced Alexander
polynomial of an oriented link $l$ is the polynomial of a
single variable obtained from 
$\D_l$ by setting all of the variables equal.) A proof can be
found in [{\bf Hir98}] (see also  [{\bf GH99}]). \bs

\ni {\bf Example 5.8.} The torus knot $5_1$ is
equivalent to $l(-2,3,1)$. Consider the
$2$-component link $l$ obtained from it by encircling the 
third band of the pretzel knot, as in Figure 8. Its Alexander
polynomial $\D_l$ is
$u_1^2-u_1^3+u_1^5 + u_2 (1-u_1^2+u_1^3)$. Using Lemma 4.1,
one easily sees that $\D_l$ has the same Mahler measure as
$u^3-u-1$, namely $\theta_0 (\approx 1.32471)$. 
\bs
\epsfxsize=1.9truein  
\centerline{\epsfbox{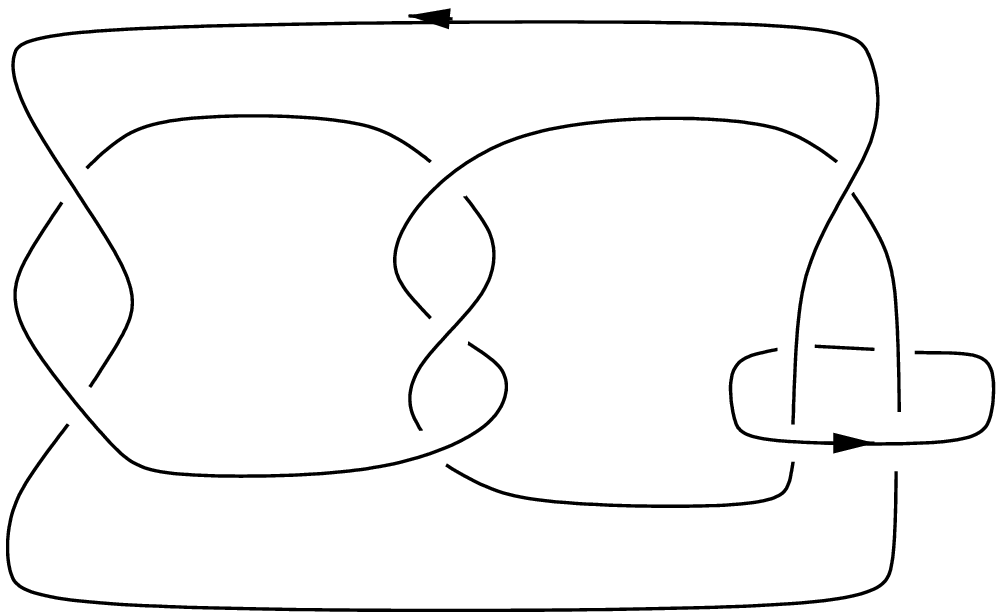}}\bs

\centerline{{\bf Figure 8:} Encircled pretzel link}\bs

We twist the two encircled arcs of $k$, forming  knots $l(q)$
as in section 2. The Mahler measures $M(\D_{l(q)})$, which
converge to $\theta_0$ by Theorem 2.2, must all be Salem
numbers in view of [{\bf Hir98}, Proposition 3.1]. 
Some of the values are given below. \bs

\centerline{\vbox{\offinterlineskip
\hrule
\halign{&\vrule#&
   \strut\quad\hfil#\quad\cr
height2pt&\omit&&\omit&\cr
&Number of Twists\hfil&&Salem Number&\cr
height2pt&\omit&&\omit&\cr
\noalign{\hrule}
height2pt&\omit&&\omit&\cr
&1\hfil &&$1\hfil $ &\cr
&2\hfil &&$1\hfil$ &\cr
&3\hfil &&$1.17628\ldots$ &\cr
&4\hfil &&$1.26123\ldots$ &\cr
&5\hfil &&$1.29348\ldots$ &\cr
&6\hfil &&$1.30840\ldots$ &\cr
&7\hfil &&$1.31591\ldots$ &\cr
&8\hfil &&$1.31986\ldots$ &\cr
&9\hfil &&$1.32201\ldots$ &\cr
&10\hfil &&$1.32319\ldots$ &\cr
&11\hfil &&$1.32385\ldots$ &\cr
&12\hfil &&$1.32423\ldots$ &\cr
height2pt&\omit&&\omit&\cr}
\hrule}} \bs

As we remarked in section 2, Theorem 2.2
bears a resemblance of form to a theorem of Thurston about
volumes of hyperbolic manifolds. In their census of the simplest
hyperbolic knots [{\bf CDW98}], P. Callahan, J. Dean and J.
Weeks record the fact that the pretzel knot $l(-2,3,7)$ has a
complement composed of only 3 ideal tetrahedra, thereby
qualifying it for the  honor of second simplest hyperbolic knot,
after the figure eight knot $4_1$ and  alongside the knot
$5_2$. \bs

\ni {\bf Question 5.9.} Is there a significant relationship
between the Mahler measure of the Alexander polynomial of
a hyperbolic link and the hyperbolic volume of its exterior?\bs

\ni {\bf Acknowledgements.} We are endebted to David Boyd, John
Dean, Michael Mossinghoff, Doug Lind, Chris
Smyth and Jeffrey Vaaler for valuable discussions, and to
Boyd, Dean and Mossinghoff for generously sharing with us the
results of computer calculations. We are grateful to the
University of Maryland for its hospitality. Finally, we
thank the referee for suggestions that improved the clarity
of our presentation.
\bs


\centerline{\bf References.} 
\bs
\baselineskip=12 pt
\item{[{\bf Ada00}]} C.~C. Adams, ``Cusp volumes and cusp densities
for hyperbolic $3$-manifolds,'' talk at 959th Amer.\ Math.\ Soc.
meeting, Columbia Univ., November 4--5, 2000.
\ss
\item{[{\bf Ahl66}]} L.~V. Ahlfors, Complex analysis, 2nd
edition, McGraw-Hill, New York, 1966.
\ss
\item{[{\bf Boy78}]} D.~W. Boyd, ``Variations on a theme
of Kronecker,'' {\sl Canad.\ Math.\ Bull.\ \bf21} (1978),
129--133.
\ss
\item{[{\bf Boy81}]} D.~W. Boyd, ``Speculations concerning the
range of Mahler's measure,'' {\sl Canad.\ Math.\ Bull.\ \bf
24} (1981), 453--469.
\ss
\item{[{\bf Boy81$'$}]} D.~W. Boyd, ``Kronecker's theorem and
Lehmer's problem for polynomials in several variables,''
{\sl J.\ Number\ Theory\ \bf 13} (1981), 116--121.
\ss
\item{[{\bf Boy98}]} D.~W. Boyd, ``Uniform approximation to Mahler's
measure in several variables,'' {\sl Canad.\ Math.\ Bull.\ \bf41}
(1998), 125--128.
\ss
\item{[{\bf Boy00}]} D.~W. Boyd, private correspondence.
\ss
\item{[{\bf CDW98}]} P.~J. Callahan, J.~C. Dean and J.R. Weeks,
``The simplest hyperbolic knots,'' preprint, 1998.
\ss
\item{[{\bf CF63}]} R.~H. Crowell and R.~H. Fox, An Introduction 
to Knot Theory, Ginn and Co., 1963.
\ss
\item{[{\bf Den97}]} C. Deninger, ``Deligne periods of mixed
motives, K-theory and the entropy of certain $Z^n$-actions,''
{\sl J.\ Amer.\ Math.\ Soc.\ \bf 10} (1997), 259--281.
\ss
\item{[{\bf Dun00}]} N.~M. Dunfield, ``The virtual Haken
conjecture: Experiments and Examples,'' talk at Barnard College
Workshop on  Topology of $3$-Manifolds, November 3, 2000.
\ss
\item{[{\bf EW99}]} G. Everest, T. Ward, Heights of
polynomials and entropy in algebraic dynamics, Springer-Verlag,
London, 1999. 
\ss
\item{[{\bf Fox60}]} R.~H. Fox, ``Free differential calculus, V,''
{Annals\ Math.\ \bf 71} (1960), 408--422.
\ss
\item{[{\bf GH99}]} E. Ghate and E. Hironaka, ``The arithmetic
and geometry of Salem numbers,'' {\sl Bull.\  Amer.\
Math.\ Soc.}, to appear.
\ss
\item{[{\bf Hir98}]} E. Hironaka, ``The Lehmer polynomial and
pretzel knots,'' {\sl Bull.\ Canad.\ Math.\ Soc.\ \bf38}
(2001), 293--314.
\ss
\item{[{\bf HW89}]} M. Hildebrand and J. Weeks, ``A computer
generated cenus of cusped hyperbolic $3$-manifolds,'' Computers
and Mathematics, E. Kaltofen and S. Watts, editors,
Springer-Verlag, Berlin, New York, 1989, 53--59.
\ss
\item{[{\bf Kaw96}]} A. Kawauchi, A Survey of Knot Theory,
Birkh\"auser, Basel, 1996. 
\ss
\item{[{\bf Kid82}]} M.~E. Kidwell, ``Relations between the
Alexander polynomial and summit power of a closed braid,''
{\sl Math.\ Sem.\ Notes\ Kobe\ Univ.\ \bf 10} (1982), 387--409.
\ss 
\item{[{\bf Kir97}]} R. Kirby, ``Problems in low-dimensional
topology.'' In Geometric Topology, W.~H. Kazez, editor, Studies
in Advanced Mathematics, A.M.S., 1997. 
\ss
\item{[{\bf Law83}]} W.~M. Lawton, ``A problem of Boyd
concerning geometric means of polynomials,'' {\sl J.\ Number\
Theory\ \bf 16} (1983), 356--362.
\ss
\item{[{\bf Leh33}]} D.~H. Lehmer, ``Factorization of certain
cyclotomic functions,'' {\sl Annals\ of\ Math.\ \bf 34} (1933),
461--479. 
\ss
\item{[{\bf Lev67}]} J. Levine, ``A method for generating link
polynomials,'' {\sl Amer.\ J.\ Math.\ \bf 89} (1967), 69--84.
\ss
\item{[{\bf Mah60}]} K. Mahler, ``An application of Jensen's
formula to polynomials,'' {\sl Mathematika\ 
\bf 7} (1960), 98--100.
\ss
\item{[{\bf Mah62}]} K. Mahler, ``On some inequalities for
polynomials in several variables,''  {\sl J.\ London\ Math.\
Soc.\ \bf 37}(1962), 341--344.
\ss
\item{[{\bf Mos98}]} M.~J. Mossinghoff, ``Polynomials with
small Mahler measure,'' {\sl Math.\ Comp.\ \bf 67} (1998),
1697--1705. 
\ss
\item{[{\bf NZ85}]} W.~D. Neumann and D. Zaiger, ``Volumes of
hyperbolic three-manifolds,'' {\sl Topology\ \bf 24} (1985), 
307--322.
\ss
\item{[{\bf Rol76}]} D. Rolfsen, Knots and Links, Publish or
Perish, Berkeley, CA, 1976. 
\ss
\item{[{\bf Sal44}]} R. Salem, ``A remarkable class of algebraic
integers. Proof of a conjecture of Vijayaraghavan,''
{\sl Duke\ Math.\ Jour.\ \bf 11} (1944), 103 -- 108. \ss
\item{[{\bf Sch95}]} K. Schmidt, Dynamical Systems of Algebraic
Origin, Birkh\"auser Verlag, Basel, 1995.
\ss
\item{[{\bf Sie44}]} C.~L. Siegel, ``Algebraic integers whose
conjugates lie in the unit circle,'' {\sl Duke\ Math.\ Jour.\ \bf
11} (1944), 597--602.
\ss
\item{[{\bf SW99}]} D.~S. Silver and S.~G. Williams, ``Mahler
measure, links and homology growth,'' {\sl Topology\ \bf41}
(2002), 979--991.
\ss
\item{[{\bf Smy00}]} C.~J. Smyth, preprint, 2000.
\ss
\item{[{\bf SM82}]} C.~J. Smyth,
``On measures of polynomials in several variables,'' {\sl Bull.\
Austral.\ Math.\ Soc.\ \bf 23} (1981), 49--63;
Corrigendum by Smyth and G. Myerson, {\sl Bull.\ Austral.\
Math.\ Soc.\ \bf 26} (1982), 317--319.
\ss
\item{[{\bf Ste78}]} C.~L. Stewart, On a theorem of Kronecker and
a related question of Lehmer, S\'eminare de Th\'eorie
de Nombres Bordeaux 1977/78, Birkha\"user, Basel, 1978.
\ss
\item{[{\bf Thu83}]} W.~P. Thurston, The Geometry and Topology
of $3$-Manifolds, Lecture Notes, Princeton 1977. Revised version
Princeton University Press, 1983.
\ss
\item{[{\bf Wal80}]} M. Waldschmidt, ``Sur le produit des
conjugu\'es ext\'erieurs au cercle unit\'e d'un entier
alg\'ebrique,'' {\sl L'Enseign.\ Math.\ \bf 26} (1980), 201--209.
\ss
 
\bs
\item{} Dept. of Mathematics and Statistics, Univ. of South Alabama, Mobile, AL  36688-0002
e-mail: silver@jaguar1.usouthal.edu,
swilliam@jaguar1.usouthal.edu

\end